\documentclass{elsart}

\usepackage{amsmath,amssymb}
\usepackage{graphicx}
\usepackage{epsfig}
\usepackage{color}
\usepackage{ulem}

\newcommand{\added}[1]{\textcolor{red}{#1}}
\newcommand{\cqfd}%
{\mbox{}\nolinebreak\hfill\rule{2mm}{2mm}\medskip\par}

\newcommand{\prodscal}[2]{\left( \, {#1} \, | \, {#2} \, \right)}

\newcommand{\derp}[2]{\frac{\partial{#1}}{\partial{#2}}}

\def\cqfd{\qed}
\def\IR{\mathbb{R}}

\def\divg{\mathop{\rm div}\nolimits}
\def\vtld{\widetilde{v}}

\def\nabtld{\widetilde{\nabla}}

\def\divtld{\widetilde{\divg}}

\newtheorem{theo}{Theorem}[section]
\newtheorem{lemm}[theo]{Lemma}

\theoremstyle{definition}

\theoremstyle{remark}

\newenvironment{dem}[1] {\par\noindent{\it Proof. }{#1}}{$\square$}
\begin{document}

\begin{frontmatter}



\title{Computing the time-continuous Optimal Mass Transport Problem
 without Lagrangian techniques}


\author[lab1]{Olivier Besson}
\ead{Olivier.besson@unine.ch}
\author[lab2]{Martine Picq}
\ead{martine.picq@insa-lyon.fr}
\author[lab2]{J\'er\^ome Pousin\corauthref{cor1}}
\corauth[cor1]{Corresponding Author Fax: 00 33 4 72 43 85 29}
\ead{jerome.pousin@insa-lyon.fr}

\address[lab1]{ Universit\'e de Neuch\^atel, Institut de Math\'ematiques
\\11, rue E. Argand, 2009 Neuch\^atel, Switzerland,}
\address[lab2]{Universit\'e de Lyon CNRS  \\INSA-Lyon ICJ UMR 5208,
bat. L. de Vinci,
\\20 Av. A. Einstein, F-69100 Villeurbanne Cedex France }

\begin{abstract}
This work originates from a heart's images tracking which is to
generate an apparent continuous motion, observable through intensity
variation from one starting image to an ending one both supposed
segmented. Given two images $\rho_0$ and $\rho_1$, we calculate an
evolution process $\rho(t,\cdot)$ which transports $\rho_0$ to
$\rho_1$ by using the optimal extended optical flow. In this paper we propose
an algorithm based on a fixed point formulation and a time-space least squares
formulation of the mass conservation equation for computing the optimal mass
transport problem. The strategy is implemented
in a 2D case and numerical results are presented with a first order Lagrange
finite element, showing the efficiency of the proposed strategy.
\end{abstract}

\end{frontmatter}

\section{Introduction}
\label{Intro} \noindent Modern medical imaging modalities can
provide a great amount of information to study the human anatomy and
physiological functions in both space and time. In cardiac Magnetic
Resonance Imaging (MRI) for example, several slices can be acquired
to cover the heart in 3D and at a collection of discrete time
samples over the cardiac cycle.  From these partial observations, the challenge is to
extract the heart's dynamics from these input spatio-temporal data throughout the cardiac cycle
\cite{Lynch08}, \cite{Schaerer08}.

\noindent Image registration consists in estimating a transformation which insures the warping
of one reference image onto another target image (supposed to present some similarity).
Continuous transformations are privileged, the sequence of transformations during the
estimation process is usually not much considered. Most important is the final resulting
 transformation and not the way one image will be transformed to the other.
Here, we consider a reasonable registration process to continuously map the image intensity
functions between two
images in the context of cardiac motion estimation and modeling.

\noindent The aim of this paper is to present, in the context of
extended optical flow, an algorithm to compute the optimal time dependent transportation
plan without using Lagrangian techniques.

\noindent The paper is organized as follows. The introduction is ended, by recalling
the optimal extended optical flow model (OEOF) . In section 2, the algorithm we propose is
presented.
Its convergence is discussed. In section 3 it is proved that solutions obtained with the
proposed algorithm are solutions to the optimal extended optical flow, that is to say to
the time dependent optimal mass transportation problem.  Section 4 deals with
numerical results. A 2D cardiac  medical image  is considered.


\subsection{The OEOF method}
\noindent 
Let us denote by $\rho$  the intensity function, and by  $v$ the
velocity of the apparent motion of brightness pattern. An image
sequence is considered via the gray-value map $\rho : Q=(0,1) \times
\Omega \rightarrow \IR$  where $\Omega \subset \IR^d$ is a bounded
regular domain, the support of images, for $d=1, 2, 3$. If image
points move according to the velocity field $v: \, Q \rightarrow
\IR^d$, then gray values $\rho(t,X(t,x))$ are constant along motion
trajectories $X(t,x)$. One obtains the optical flow equation:
\begin{equation}
 \frac{d}{dt} \rho(t,X(t,x))= \partial_t \rho(t,X(t,x))  +
 \prodscal{v}{\nabla_X \rho(t,X(t,x))}_{\IR^d}  = 0.
\end{equation}

The assumption that the pixel intensity does not change during the
movement is in some cases too restrictive. A weakened assumption
sometimes called extended optical flow, can replace the intensity
preservation
 by a mass preservation condition which reads:
\begin{equation}
\partial_t \rho + \prodscal{v}{\nabla_x \rho}_{\IR^d} + \divg{(v)} \rho =
0.
\end{equation}
The previous equations lead to an ill-posed problem for the unknown
$(\rho,v)$. Variational formulations or relaxed minimizing problems
for computing jointly $(\rho,v)$ have been first proposed in
\cite{Aubert:1999} and after by many other authors. Here our concern
is somewhat  different. Finding $(\rho,v)$ simultaneously  is
possible by solving the optimal mass transport problem
\eqref{optiproblem}-\eqref{optiproblem2}, developed in
\cite{Benamou-Brenier,Brenier}.

\noindent Let $\rho_0$ and $\rho_1$ be the cardiac images
between two times arbitrary fixed to zero and one, the mathematical
problem reads: find $\rho$ the gray level function defined from $Q$
with values in $[0,1]$ verifying
\begin{equation}\label{optiproblem}
\left \{ \begin{array}{l}
\partial_t \rho(t,x)  + \divg(v(t,x)\rho(t,x)) = 0, \, \text{ in }
(0,1)\times \Omega\\
\rho(0,x)=\rho_0(x); \quad \rho(1,x)=\rho_1(x)\\
\end{array} \right.
\end{equation}
 The
velocity function $v$ is determined in order to minimize the
functional:
\begin{equation}\label{optiproblem2}
\inf_{\rho,v} \int_0^1 \int_\Omega \rho(t,x) \| v(t,x) \|^2 \, dtdx.
\end{equation}
Thus we get an image sequence through the gray-value map $\rho$.
Let us mention \cite{Tannenbaum}, for example,  where the optimal
mass transportation approach is used in images processing. For general
properties of optimal
 transportation, the reader is referred to the books by C. Villani \cite{Vil07}
and L. Ambrosio et al. \cite{Ambrosio2}.

\section{Algorithm for solving the Optimal Extended Optical Flow}

\noindent In what follows, let us specify our hypotheses.
\begin{itemize}
\item [H1] $\Omega$ is a bounded $C^{2,\alpha}$ domain satisfying the
exterior sphere condition.

\item [H2]  $\rho_i \in C^{1,\alpha}(\overline \Omega)$ for $i=1,2$,  and
$\rho_0 = \rho_1 \text{ on }\partial \Omega$.
Moreover there exist two constants such that
$0<\underline \beta \le \rho_i \le \overline \beta$ in $\Omega$.
\end{itemize}

\noindent Let $\rho^0 \in C^{1,\alpha}([0,1]\times\overline \Omega)$
be given by $\rho^0(t,x)=(1-t)\rho_0(x) + t \rho_1(x)$. We have
$\|\partial_t \rho^0\|_{C^{0,\alpha}([0,1]\times\overline \Omega)}\le C(\rho_0, \rho_1)$
and $\partial_t \rho^0\vert_{\partial \Omega}=0$.

\noindent For each $t\in [0,1]$, our need for problem \eqref{optiproblem}-\eqref{optiproblem2} is a
velocity field vanishing on $\partial \Omega$. To do so, the following method is used.
\begin{itemize}
\item
Compute
\begin{equation}\label{problemC}
\left \{ \begin{array}{l}
 -\divg(\rho^{n}(t,\cdot)\nabla \eta) = 0 \, \text{ in }
 \Omega\\
\rho^{n}(t,\cdot)\partial_n \eta=1 \quad \text{ on }
 \partial \Omega,
\end{array} \right.
\end{equation}
and set
$ C^n(t) = \frac{1}{\vert \partial \Omega \vert} \int_\Omega \partial_t \rho^{n} \eta \, dx$.

\item For each $t\in [0,1]$ compute $\varphi^{n+1}$  solution to

\begin{equation}\label{problemphi}
\left \{ \begin{array}{l}
 -\divg(\rho^{n}(t,\cdot)\nabla \varphi^{n+1}) = \partial_t \rho^{n}(t,\cdot), \, \text{ in }
 \Omega\\
\varphi^{n+1}=C^n(t) \quad \text{ on }
 \partial \Omega.
\end{array} \right.
\end{equation}

\item Set $v^{n+1}= \nabla \varphi^{n+1}$.

\item Compute $\rho^{n+1}$, \added{$L^2$-least squares} solution to
\begin{equation}\label{problemrho}
\left \{ \begin{array}{l}
\partial_t \rho^{n+1}(t,x)  + \divg(v^{n+1}(t,x)\rho^{n+1}(t,x)) = 0, \, \text{
in }
(0,1)\times \Omega\\
\rho^{n+1}(0,x)=\rho_0(x); \quad \rho^{n+1}(1,x)=\rho_1(x).\\
\end{array} \right.
\end{equation}

\end{itemize}

\noindent For each $t \in [0,1]$, since $\rho^n(t,\cdot)$, and $\partial_t \rho^n (t,\cdot)\in
C^{0,\alpha}(\overline \Omega)$, theorem 6.14 p. 107 of \cite{Gilbarg} applies, and
there exists a unique $\varphi^{n+1}(t,\cdot) \in C^{2,\alpha}(\overline
\Omega)$ solution of
problem \eqref{problemphi}. \added{In problem \eqref{problemphi} the  time is a
parameter. As the following regularities with respect to time are verified:
$\rho^n \in C^{1,\alpha}; \, \partial_t \rho^n \in C^{0,\alpha}; \,  C^n \in
C^{0,\alpha} $. The classical $C^{2,\alpha}(\overline \Omega)$ a priori
estimates for solutions to elliptic problems allow us to prove that
$\varphi^{n+1}$ is a $C^{0,\alpha} $ function with respect to time.  So we
have:}
$$
\|\varphi^{n+1}\|_{C^{0,\alpha}([0,1];C^{2,\alpha}(\overline \Omega))}\le M(\|C^n
\|_{C^{0,\alpha}([0,1])} + \|\partial_t \rho^n\|_{C^{0,\alpha}([0,1]\times\overline \Omega)}).
$$
Consider the extension of $\varphi^{n+1}$ by $C^n$ outside of the domain $\Omega$; still denoted by
$\varphi^{n+1}$. Since the right hand side of equation \eqref{problemphi} vanishes on
$\partial\Omega$, this extension is regular, and  the function $v^{n+1}$ vanish
outside $\Omega$ and belongs to $C^{0,\alpha}([0,1];C^{1,\alpha}(\IR^2))$.

Define the two flows
$X^{n+1}_{\pm}(s,t,x)\in C^{1,\alpha}([0,1]\times[0,1]\times\IR^2 ; \IR^2) $ by
\begin{equation}\label{problemflo}
\left \{ \begin{array}{l}
\frac{d}{ds}  X^{n+1}_{\pm}(s,t,x)  = \pm v^{n+1}(s,X^{n+1}_{\pm}(s,t,x))\, \text{ in }
(0,1)\\
X^{n+1}_{\pm}(t,t,x)  = x.\\
\end{array} \right.
\end{equation}

We have the following

\begin{lemm}\label{existrho}
The $L^2$-least squares solution to problem \eqref{problemrho} is given by:
\begin{equation}\label{repform}
\begin{array}{c}
\rho^{n+1}(t,x)= (1-t) \frac{\rho^2_0(X^{n+1}_{+}(0,t,x))}{\rho^{n}(t,x)}\\ +
t \frac{\rho^2_1(X^{n+1}_{+}(1,t,x))}{\rho^{n}(t,x)}.
\end{array}
\end{equation}
Moreover, if $0<\underline \beta \le \rho^{n} \le \overline \beta$ in
$[0,1]\times \overline \Omega$, then
$\rho^{n+1} \in C^{1,\alpha}(\added{ [0,1] \times \Omega)}$, and verifies the same property.
\end{lemm}
\begin{dem}
We have $X^{n+1}_{-}(1-s,1-t,x)=X^{n+1}_{+}(s,t,x)$ for every
$(s,t,x)\in [0,1]\times[0,1]\times\IR^2$ (see for example \cite{Ambrosio}).

Let us express equation \eqref{problemrho} along the integral curves of equation
\eqref{problemflo}.
The $L^2$-least squares solution to the  ordinary differential equation with
initial and final conditions reads
\begin{equation}\begin{array}{c}\label{repres}
\rho^{n+1}(s,X^{n+1}_{+}(s,t,x)))=
(1-s)e^{-\int_0^s \divg(v^{n+1}(\tau,X^{n+1}_{+}(\tau,t,x)))\, d\tau} \rho_0(X^{n+1}_{+}(0,t,x)) \\
+ s e^{\int_s^1 \divg(v^{n+1}(\tau,X^{n+1}_{+}(\tau,t,x)))\, d\tau}\rho_1(X^{n+1}_{+}(1,t,x)).
\end{array}
\end{equation}
Equation \eqref{problemphi} gives the following expression for the divergence
\begin{multline}
\divg(v^{n+1}(s,X^{n+1}_{+}(s,t,x))= \divg(v^{n+1}(s,X^{n+1}_{-}(1-s,1-t,x))\\
=\frac{d}{ds} \ln (\rho^n(s,X^{n+1}_{-}(1-s,1-t,x))).
\end{multline}
The representation formula \eqref{repform} is straightforwardly deduced from \eqref{problemphi}.
The regularity of the function $\rho^{n+1}$ is a consequence of the regularity of the flow $X_+^{n+1}$.
\end{dem}

\noindent Let us now consider the convergence of the algorithm \eqref{problemC}-\eqref{problemrho}.
\begin{theo}\label{convalgo}
There exist
$(\rho,\varphi)\in C^{1}([0,1]\times\overline \Omega)\times
C^{0}([0,1];C^{2}(\overline \Omega))$,
\added{$L^2$-least squares solution, respectively} solution to
\begin{equation}\label{problemlim1a}
\left \{ \begin{array}{l}
\partial_t \rho(t,x)  + \divg(\nabla \varphi(t,x)\rho(t,x)) = 0, \; \mathrm{in}
\,
(0,1)\times \Omega\\
\rho(0,x)=\rho_0(x); \quad \rho(1,x)=\rho_1(x) \; \mathrm{in} \, \Omega \\
\end{array} \right.
\end{equation}
\begin{equation}\label{problemlim1b}
\left \{ \begin{array}{l}
 -\divg(\rho(t,\cdot)\nabla \varphi) = \partial_t \rho(t,\cdot), \; \mathrm{in} \,
 \Omega\\
\varphi=C(t);\, \nabla  \varphi = 0 \; \mathrm{on} \,
 \partial \Omega\\
\end{array} \right.
\end{equation}
with $C(t)$ defined by:
\begin{equation}\label{problemlim2}
\left \{ \begin{array}{l}
 -\divg(\rho(t,\cdot)\nabla \eta) = 0 \; \mathrm{in} \, \Omega\\
\rho(t,\cdot) \, \partial_n \eta=1 \; \mathrm{on} \,
 \partial \Omega\\
 C = \frac{1}{\vert \partial \Omega \vert} \int_\Omega \partial_t \rho \, \eta \, dx.
\end{array} \right.
\end{equation}
\end{theo}
\begin{dem}
Since
$\|v^0 \|_{C^{0,\alpha}([0,1])} + \|\partial_t \rho^0\|_{C^{0,\alpha}([0,1]\times\overline
\Omega)}$
is bounded,
$\|\varphi^{n+1}\|_{C^{0,\alpha}([0,1];C^{2,\alpha}(\overline \Omega))}$
and $\|v^{n+1}\|_{C^{0,\alpha}([0,1];C^{1,\alpha}(\IR^2))}$ are uniformly bounded in $n$.

\noindent From lemma \ref{existrho} there exists a unique
$\rho^{n+1}$, the $L^2$-least squares solution of \eqref{problemrho}.
Let us give an estimate for $D_3X_{+}^{n+1}$.
Starting from
$$
D_1 X^{n+1}_{+}(s,t,x))= v^{n+1}(s,X^{n+1}_{+}(s,t,x)),
$$
we deduce (see \cite{Ambrosio})
\begin{equation}
\left \{ \begin{array}{l}
D_3D_1 X^{n+1}_{+}(s,t,x)  = D_2 v^{n+1}(s,X^{n+1}_{+}(s,t,x))D_3X^{n+1}_{+}(s,t,x)\\
D_3X^{n+1}_{+}(t,t,x)  = Id.\\
\end{array} \right.
\end{equation}
Since $D_3D_1X^{n+1}_{+}(s,t,x)  =D_1D_3X^{n+1}_{+}(s,t,x)$ we get
\begin{equation}\label{repDX}
D_3X^{n+1}_{+}(s,t,x)  = e^{-\int_t^s D_2(v^{n+1}(\tau,X^{n+1}_{+}(\tau,t,x)))\, d\tau} Id.
\end{equation}
Thus $\|D_3v_{+}^{n+1}\|_{C^{0,\alpha}([0,1]^2\times \IR^2)}$ is uniformly bounded in $n$.

\noindent Since we have \cite{Ambrosio}:
$$
D_2X^{n+1}_{+}(s,t,x)= \prodscal{v^{n+1}(s,t,x)}{D_3X^{n+1}_{+}(s,t,x)}
$$
we obtain a bound for $\|D_2v^{n+1}\|_{C^{0,\alpha}([0,1]^2\times \IR2)}$ independent of $n$.

\noindent From theorem \ref{existrho} we deduce that
$\|\rho^{n+1}\|_{C^{1,\alpha}([0,1]\times  \overline \Omega)}$ is uniformly bounded. Since the
embeddings
$$C^{0,\alpha}([0,1];C^{2,\alpha}(\overline \Omega))\hookrightarrow
C^{0}([0,1];C^{2}(\overline \Omega))\, \mathrm{and} \,  C^{1,\alpha}([0,1]\times
 \overline
\Omega)\hookrightarrow C^{1}([0,1]\times  \overline \Omega)
$$
are relatively compact there is a subsequence of $(\rho^{n},\varphi^{n})$
solution to \eqref{problemC}-\eqref{problemrho}, still denoted by $(\rho^{n},\varphi^{n})$
converging to
$(\rho,\varphi)$ in $C^{1}([0,1]\times  \overline \Omega)\times C^{0}([0,1];C^{2}(\overline
\Omega))$, and $(\rho,\varphi)$ is the solution of
\eqref{problemlim1a}-\eqref{problemlim2} \added{ provided the boundary
conditions to be  justified. The condition $\nabla \varphi^n\vert_{\partial
\Omega} =0$ is valid for the approximations $\varphi^n$ (since the functions can
be extended by $C^n$ outside of $\Omega$). So the convergence in
$C^{0}([0,1];C^{2}(\overline \Omega))$ yields the condition for the gradient of
limit function. For the approximations of function $\rho$, the formula given in
Lemma \ref{existrho} combined with the regularity result show that the boundary
conditions are exactly satisfied. These conditions are thus valid for the limit
function due to the convergence in $C^1$.}
\end{dem}

\added{We will show in the next section that the above least squares solution
$\rho$ is in fact a classical solution.}

\section{Interpretation of solutions to problem
\eqref{problemlim1a}-\eqref{problemlim2}}

\noindent In this section it is shown that the solution to problem
\eqref{problemlim1a}-\eqref{problemlim2} is a solution to the time dependent
optimal mass transportation problem.

\noindent From one hand, remark that \added{$\varphi$ solution to problem
\eqref{problemlim1b} satisfies:
$$
\varphi-C =\underset{\psi\in L^2((0,1);H^1_0(\Omega))}{\rm Argmin}\frac{1}{4} \int_0^1 \| \partial_t \rho  +
\divg(\rho \nabla \psi ) \|_{H^{-1}(\Omega)}^2 \, dt.
$$
Since the functions $(\rho,\varphi)$ are sufficiently regular, we have:
$$
\varphi-C =\underset{\psi\in L^2((0,1);H^1_0(\Omega)\cap H^2(\Omega))}{\rm Argmin} \frac{1}{4}\int_0^1 \| \partial_t \rho  +
\divg(\rho \nabla \psi ) \|_{L^{2}(\Omega)}^2 \, dt.
$$
From an other hand, zero is a bound from below of the functional to be minimized with respect to $(u,\psi)$:
$$
\begin{array}{l}
 0= \frac{1}{4}\int_0^1 \| \partial_t \rho  +
\divg(\rho \nabla (\varphi-C) ) \|_{L^{2}(\Omega)}^2 \, dt \le  \\
\underset{
\begin{array}{c}
\scriptstyle \{\psi\in L^2((0,1);H^1(\Omega)), \; u\in L^2((0,1);L^2(\Omega)) \\*[-3mm]
\scriptstyle \partial_t u  + \divg(-u \nabla \psi)) \in L^2((0,1);L^2(\Omega)) \\*[-3mm]
\scriptstyle \partial_t u  + \divg(u \nabla \psi) = 0 \\*[-3mm]
\scriptstyle \nabla \psi\vert_{\partial \Omega} = 0 \\*[-3mm]
\scriptstyle \psi\vert_{\partial \Omega} = C \\*[-3mm]
\scriptstyle u(0)=\rho_0;\; u(1)=\rho_1 \text{ in } \Omega \}
\end{array}
}
{\rm Min}
\frac{1}{4} \int_0^1 \| \partial_t u  +
\divg(u \nabla \psi ) \|_{L^2(\Omega)}^2 \, dt.
\end{array}
$$
}
We deduce that $(\rho,\varphi)$, solution to problem
\eqref{problemlim1a}-\eqref{problemlim2}, satisfies
\begin{equation}\label{problemint1}
(\rho,\varphi)=
\underset{
\begin{array}{c}
\scriptstyle \{\psi\in L^2((0,1);H^1(\Omega)), \; u\in L^2((0,1);L^2(\Omega)) \\*[-3mm]
\scriptstyle \partial_t u  + \divg(-u \nabla \psi)) \in L^2((0,1);L^2(\Omega)) \\*[-3mm]
\scriptstyle \partial_t u  + \divg(u \nabla \psi) = 0 \\*[-3mm]
\scriptstyle \nabla \psi\vert_{\partial \Omega} = 0 \\*[-3mm]
\scriptstyle \psi\vert_{\partial \Omega} = C \\*[-3mm]
\scriptstyle u(0)=\rho_0;\; u(1)=\rho_1 \text{ in } \Omega \}
\end{array}
}
{\rm Argmin}
\frac{1}{4} \int_0^1 \| \partial_t u  +
\divg(u \nabla \psi ) \|_{L^2(\Omega)}^2 \, dt.
\end{equation}
\begin{lemm}\label{problemH-1}
Let $(\rho,\varphi)$ be a solution to problem
\eqref{problemlim1a}-\eqref{problemlim2}. Then it
satisfies
\begin{equation}\label{problemint2}
(\rho,\varphi)=
\underset{
\begin{array}{c}
\scriptstyle \{\partial_t u  + \divg(u \nabla \psi)=0; \, \nabla
\psi\vert_{\partial \Omega}=0;
\\*[-3mm]
\scriptstyle \psi\vert_{\partial \Omega}=C; \, u(0)=\rho_0; \, u(1)=\rho_1 \text{ in } \Omega \}
\end{array}
}
{\rm Argmin}
\ \int_0^1 \| \divg(u \nabla \psi ) \|_{H^{-1}(\Omega)}^2 \, dt.
\end{equation}
\end{lemm}
\begin{dem} This is a simple consequence of
$\partial_t \rho=-\divg(\rho \nabla \varphi)$,
and of the regularity of $\divg(\rho \nabla \varphi)$ which implies
$\| \divg(\rho \nabla \varphi)\|_{L^2(\Omega)}=\| \divg(\rho \nabla
\varphi)\|_{H^{-1}(\Omega)}$.
\end{dem}

\begin{theo}\label{problemwas}
Let $(\rho,\varphi)$ be solution to problem
 \eqref{problemlim1a}-\eqref{problemlim2}, the existence of which is given in
 Theorem \ref{convalgo}, then it satisfies:
\begin{equation}\label{eqwas}
(\rho,\nabla\varphi)=
\underset{\{ \partial_t u  + \divg(u v)=0; \, u(0)=\rho_0; \, u(1)=\rho_1 \text{ in } \Omega\}}
{\rm Argmin}
\int_0^1 \int_\Omega u \| v \|^2 \, dxdt.\\
\end{equation}
\end{theo}
\begin{dem}
Choose $u$ regular verifying
$0< \underline \beta \le u \le \overline \beta$, and for all $t\in (0,1)$ solve
\begin{equation}\label{probleminfv}
\inf_{\{ v\in L^2(\Omega) \,\partial_t u  + \divg(u v)=0\}}   \int_\Omega u \| v \|^2 \, dx.
\end{equation}
Let $H=H^1_0(\Omega)$
be equipped with the following inner product:
$$
(\theta,\psi)=\int_\Omega u \prodscal{\nabla \theta}{\nabla \psi} \, dx,
$$
which induces a semi-norm which is equivalent to the $H^1$-norm.
The Riez's theorem claims that for the linear continuous form
$$
\mathcal{L}_u(\psi) =<-\divg{(uv)},\psi>_{H;H'}=<\partial_tu,\psi>_{H;H'},
$$
there is a unique $\theta \in H$ such that
$$
\mathcal{L}_u(\psi)=\int_\Omega u \prodscal{\nabla \theta}{\nabla \psi} \, dx,
\, \forall \psi \in H.
$$
Therefore $v=\nabla \theta$ and problem \eqref{probleminfv} is reduced to
\begin{equation}\label{probleminfv2}
\inf_{\{ \psi \in H, \, \partial_t u  + \divg(u \nabla \psi)=0, \,
\psi\vert_{\partial \Omega}=C\}}   \int_\Omega u \| \nabla \psi \|^2 \, dx.
\end{equation}

Since
$$
\int_\Omega u \| \nabla \psi \|^2 \, dx = \|\divg(u\nabla \psi)\|^2_{H'},
$$
problem \eqref{probleminfv2} reads
\begin{equation}\label{probleminfv3}
\inf_{\{ \psi \in H, \, \partial_t u  + \divg(u \nabla \psi)=0, \,
\psi\vert_{\partial \Omega}=C\}}  \|\divg(u\nabla \psi)\|^2_{H'}
\end{equation}
or
\begin{equation}\label{probleminfv4}
\inf_{\{ \psi \in H, \, \partial_t u  + \divg(u \nabla \psi)=0, \,
\psi\vert_{\partial \Omega}=C\}} \frac{1}{4} \| \partial_tu + \divg(u\nabla
\psi)\|^2_{H'}.
\end{equation}
Gathering lemma \ref{problemH-1} with the previous result proves the theorem.
\end{dem}

\section{Numerical Approximation of the 2D Optimal Extended Optical Flow}

\noindent The numerical method is based on a finite element time-space
$L^2$ least squares formulation (see \cite{Besson}) of the linear conservation law
\eqref{problemrho}.

\noindent Define $\vtld^{n+1}$ as
\[\vtld^{n+1} = (1,v^{n+1}_{1},v^{n+1}_{2})^t \]
and for  a sufficiently regular function $\varphi$ defined on $Q$,
set
\[\nabtld\varphi =
\left( \derp{\varphi}{t}, \derp{\varphi}{x_{1}},
\derp{\varphi}{x_{2}} \right)^t,
\]
and
\[
\divtld(\vtld^{n+1}  \ \varphi) = \derp{\varphi}{t} +  \sum_{i=1}^{2}
\derp{ }{x_{i}}( v^{n+1}_{i} \ \varphi).
\]
Let $\{\varphi_1 \cdot \cdot \cdot \varphi_N\}$ be a basis of a space-time finite element
subspace
\[
V_h = \{\varphi, \text{ piecewise regular polynomial functions, with } \varphi(0,\cdot)=
\varphi(1,\cdot)=0 \},
\]
for example, a  brick Lagrange finite element of order one (\cite{Besson2}).
Let $\Pi_h$ be the Lagrange interpolation operator. Let also $W_h$ be the
finite element subspace of $H^1_0(\Omega)$, where the basis functions
$\{\psi_1 \cdot \cdot \cdot \psi_M\}$ are  the traces at $t=0$ of basis functions
$\{\varphi_i\}_{i=1}^N$.
An approximation of problem \eqref{problemphi} is: for a discrete sequence of time $t$ compute
\begin{equation}
\int_\Omega(\rho_h^{n}(t,\cdot)\prodscal{\nabla (\varphi_h^{n+1}-C^n(t))}{\nabla
\psi_h}\, dx =
\int_\Omega \partial_t \rho_h^{n}(t,\cdot) \psi_h \, dx \quad \forall  \psi_h \in W_h,
\end{equation}
and define  $\vtld^{n+1}= \nabla \varphi_h^{n+1}$.
The $L^2$ least squares formulation of problem \eqref{problemrho} is defined in the following
way. Consider the functional
\[
J(c) = \frac{1}{2}  \int_{Q} \left(\divtld(\vtld^{n+1}  \ c)
+\partial_t \rho^{n}_h + \divtld\left[\vtld^{n+1} \ \Pi_h \big ((1-t)\rho_{0} +
t\rho_{1}\big )\right]\right )^2 \, dx \, dt.
\]
This functional is convex and coercive in an appropriate anisotropic
Sobolev's space \cite{Besson}. The minimizer  of $J$ is
$\rho^{n+1}_h-\Pi_h \big ((1-t)\rho_{0} + t\rho_{1}\big )$ which is the solution
to the  following problem
\begin{multline}\label{dirihomoh}
\int_{Q} \divtld(\vtld^{n+1} \ \rho^{n+1}_h) \cdot \divtld(\vtld^{n+1} \ \psi_h) \, dx \
dt = \\
\int_{Q} \left(-\partial_t \rho_h^n -  \divtld\left(\vtld^{n+1} \, \Pi_h \big ((1-t)\rho_{0} +
t\rho_{1}\big )\right)\right) \cdot \divtld(\vtld^{n+1} \ \psi_h) \, dx \ dt
\end{multline}
for all $\psi_h  \in V_h$, where
$$
\rho_h = \sum_{i=1}^N  \rho_i \varphi_i(t,x).
$$
Thus an approximation of the solution to
problem \eqref{problemrho} is $\rho^{n+1}_h - \Pi_h \big ((1-t)\rho_{0} +
t\rho_{1}\big )\in V_h$.

\noindent The iterative strategy  described in Section 2 is used to compute an approximated solution, and to
reconstruct the systole to diastole images of a slice of a left ventricle.
Ten time steps have been used to compute the solution, and 10000 degrees
of freedom for the time-space least squares finite element. The approximated
fixed point algorithm converges in about 10 iterations with an accuracy of about
$10^{-7}$. 
In the next figure \ref{systole}, the initial image and the final image are presented.
\begin{figure}[h]
\begin{center}
\begin{minipage}{0.48\linewidth}
\epsfig{file=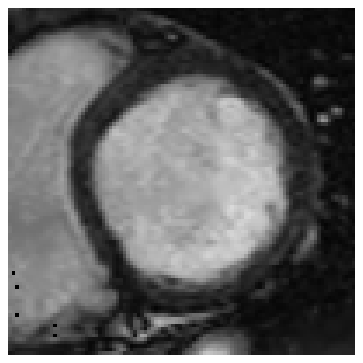, width=5cm}
\end{minipage}
\begin{minipage}{0.48\linewidth}
\epsfig{file=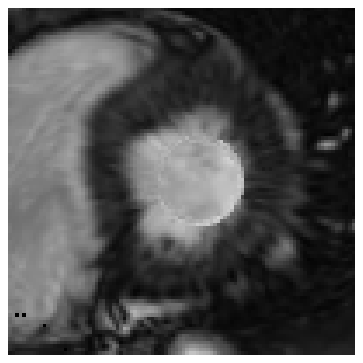,width=5cm}
\end{minipage}
\caption{End of diastole of a left ventricular (a), of systole (b) }\label{systole}
\end{center}
\end{figure}
In the following  figure \ref{syst_5}, two intermediate times $1/3$ and $2/3$ are shown.
\begin{figure}[h]
\begin{center}
\begin{minipage}{0.48\linewidth}
\epsfig{file=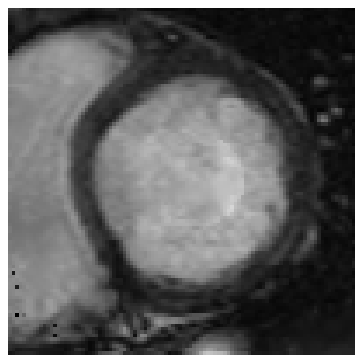, width=5cm}
\end{minipage}
\begin{minipage}{0.48\linewidth}
\epsfig{file=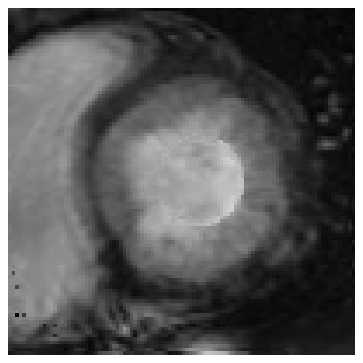,width=5cm}
\end{minipage}
\caption{Time step 3 and 6}\label{syst_5}
\end{center}
\end{figure}

\noindent To summarize, in this work, we present a fixed point algorithm
for the computation of the time dependent optimal mass transportation problem, allowing to  handle
the images tracking
problem. The efficiency of the method has been tested with a 2D example.

\end{document}